\documentclass{elsart}
\usepackage{graphics}
\usepackage{graphicx}
\usepackage{epsfig}
\usepackage{amssymb,amsmath}
\usepackage[colorlinks]{hyperref}

\newtheorem{lemma}{Lemma}[section]

\newtheorem{theorem}{Theorem}[section]
\newtheorem{proposition}{Proposition}[section]
\newtheorem{definition}{Definition}[section]

\begin{document}
\begin{frontmatter}
\title{Some Generalized Clifford-Jacobi Polynomials and Associated Spheroidal Wavelets}
\author{Sabrine Arfaoui}
\address{Department of Informatics, Higher Institute of Applied Sciences and Technology of Mateur, Street of Tabarka, 7030 Mateur, Tunisia.\\
and\\
Research Unit of Algebra, Number Theory and Nonlinear Analysis UR11ES50, Faculty of Sciences, Monastir 5000, Tunisia.}
\ead{arfaoui.sabrine@issatm.rnu.tn}
\author{Anouar Ben Mabrouk}
\address{Higher Institute of Applied Mathematics and Informatics, University of Kairouan, Street of Assad Ibn Alfourat, Kairouan 3100, Tunisia.\\
and\\
Research Unit of Algebra, Number Theory and Nonlinear Analysis UR11ES50, Faculty of Sciences, Monastir 5000, Tunisia.}
\ead{anouar.benmabrouk@fsm.rnu.tn}
\begin{abstract}
In the present paper, by extending some fractional calculus to the framework of Cliffors analysis, new classes of wavelet functions are presented. Firstly, some classes of monogenic polynomials are provided based on 2-parameters weight functions which extend the classical Jacobi ones in the context of Clifford analysis. The discovered polynomial sets are next applied to introduce new wavelet functions. Reconstruction formula as well as Fourier-Plancherel rules have been proved. The main tool reposes on the extension of fractional derivatives, fractional integrals and fractional Fourier transforms to Clifford analysis.
\end{abstract}
\end{frontmatter}
\section{Introduction}
Clifford Algebra is characterised by additional concepts as it provides a simpler model of mathematical objects compared to vector algebra. It permits a simplification in the notations of mathematical expressions such as plane and volume segments in two, three and higher dimensions by using a coordinate-free representation. Such representation is characterised by an important feature resumed in the fact that the motion of an object may be described with respect to a coordinate frame defined on the object itself. This means that it permits to use a self-coordinate system related to the object in hand.

In the present work, one aim is to provide a rigourous development of wavelets adapted to the sphere based on Clifford calculus. The frame is somehow natural as wavelets are characterized by scale invariance of approximation spaces. Clifford algebra is one mathematical object that owns this characteristic. Recall that multiplication of real numbers scales their magnitudes according to their position in or out from the origin. However, multiplication of the imaginary part of a complex number performs a rotation, it is a multiplication that goes round and round instead of in and out. So, a multiplication of spherical elements by each other results in an element of the sphere. Again, repeated multiplication of the imaginary part results in orthogonal components. Thus, we need a coordinates system that results always in the object, a concept that we will see again and again in the Algebra. In other words, Clifford algebra generalizes to higher dimensions by the same exact principles applied at lower dimensions, by providing an algebraic entity for scalars, vectors, bivectors, trivectors, and there is no limit to the number of dimensions it can be extended to. More details on Clifford analysis, clifford calculus, origins, history, developments may be found in \cite{Abreuetal}, \cite{Delanghe}, \cite{DeBie}, \cite{Pena}, \cite{Lehar}, .

In the present work, we propose to develop new wavelet analysis constructed in the framework of Clifford analysis by adopting monogenic functions which may be described as solutions of the Dirac operator and are direct higher dimensional generalizations of holomorphic functions in the complex plane. We apply such extension to some well adapted Clifford weights to construct new spheroidal wavelets. Recall that wavelets are widespread in the last decades. They have become an interesting and useful tool in many fields such as mathematics, quantum physics, electrical engineering and seismic geology and they have proved to meet a need in signal processing that Fourier transform was not the best answer. Classical Fourier analysis provides a global description of signals and did not provide a time localization.

Wavelet analysis starts by convoluting the analyzed function with copies $\psi_{a,b}$, $a>0$ and $b\in\mathbb{R}$ (called wavelets) issued from a source (mother wavelet) $\psi$ by dilation $a>0$ and translation $b\in\mathbb{R}$, \begin{equation}\label{psiab}
\psi_{a,b}(x)=\dfrac{1}{\sqrt{a}}\psi(\dfrac{x-b}{a}).
\end{equation}
Generally, the source $\psi$ has to satisfy the so-called admissibility condition
\begin{equation}\label{admissibility-conditionofpsi}
\mathcal{A}_\psi=\displaystyle\int_{-\infty}^{+\infty}\dfrac{|\widehat{\psi}(u)|^2}{|u|}du<+\infty,
\end{equation}
where $\widehat{\psi}$ is the Fourier transform of $\psi$.

The wavelet analysis starts by computing the Continuous Wavelet Transform (CWT) of the analyzed function $f$
\begin{equation}\label{waveletcoefficientcab(f)}
C_{a,b}(f)=<f,\psi_{a,b}>=\dfrac{1}{\sqrt{a}}\,\displaystyle\int_{-\infty}^{+\infty}f(x)\overline{\psi(\dfrac{x-b}{a})}dx
\end{equation}
which in turns permits to re-construct in some sense the analyzed function $f$ via an inverse transform based on the admissibility assumption (\ref{admissibility-conditionofpsi}) as
\begin{equation}
f(x)=\dfrac{1}{\mathcal{A}_\psi}\displaystyle\int_{\mathbb{R}}\displaystyle\int_{0}^{+\infty}C_{a,b}(f)\psi_{a,b}(x)\dfrac{da}{a^2}db.
\end{equation}
In the present paper we propose to construct some special spheroidal wavelets in the context of Clifford analysis with the help of fractional calculus. Recall that fractional calculus is a generalization of ordinary differentiation and integration to arbitrary (non-integer) orders. Such topic is not new but recently a coming back to its application has taken place in various areas of engineering, science, finance, applied mathematics and bio-engineering.

The wavelets constructed here are general copies of the Gegenbauer-Jacobi one developed in \cite{Brackx-Schepper-Sommen1}, \cite{Brackx-Schepper-Sommen2}, \cite{Brackx-Schepper-Sommen3}, \cite{Brackx-Schepper-Sommen4} and \cite{Brackx-Schepper-Sommen5} where the authors tried to point out a wavelet analysis in homogenous euclidean spaces. The main idea applied there was by considering a clifford Heaviside function to decompose the weight $(1-x)^\alpha (1+x)^\beta$ into semi-radial ones and thus to apply then the radial bases of euclidiean spaces by assuming that $\alpha\pm\beta=1$. The basic idea applied there is inspired from ... and resumed in the fact that for this case the weight function may be decomposed into a sum of tow parts: one part is the well known radial classical weight $(1+|x|^2)^\alpha$ and a second semi-radial one which leads with the well known rule of derivation of integrals with parameters to the first part. This was the main difference and motivation with our case where such a decomposition is not possible.

This was the main motivation and the main difference with the present paper.+ We instead come back to fractional calculus to overcome the problem of the non radially symmetric weight applied here and did not restrict to the previous case. We will prove instead that fractional calculus may be a good tool to overcome the difficulties crossed in the new context.

The organization of this paper is as follows: In section 2 a revision of fractional calculus such as fractional derivation and fractional Fourier transform is provided. Section 3 is devoted to Clifford calculus. Basic operations, Clifford fractional derivation, Clifford Fourier transform and Fractional Clifford Fourier transform are reviewed. Next, our idea of the generalization of Clifford-Jacobi polynomials is developed in section 4. Section 5 is concerned with the development of new wavelets in the Clifford context associated to the polynomial class developed previously. Relative continuous wavelet is also provided and reconstruction rule is proved in the new framework. We concluded afterward
\section{Fractional calculus}
We propose in this part to review basic definitions of fractional differentiation as well as fractional Fourier transforms to be applied next. For backgrounds on this part, the readers may refer to \cite{Baleanu}, \cite{Das}, \cite{Duboisetal}, \cite{Janin}, \cite{Katugampola}, \cite{Kilbas}, \cite{Ortigueira}, \cite{Ross}, \cite{Sabatieretal}, \cite{Samko}, \cite{Scalesetal}, \cite{Tarasovetal}, \cite{Uchaikin}, \cite{Wang1}, \cite{Wang2}, \cite{Yang}.
\subsection{Fractional derivation}
There are in fact many ideas to introduce the fractional derivative(s) of functions. One well known method reposes on Riemman-Liouville fractional integral which consists in a natural extension of the Cauchy formula given by
\begin{equation}\label{cauchyformulaforf}
J^nf(t)=f_n(t)=\dfrac{1}{(n-1)!} \displaystyle\int_{0}^{t} (t-\tau)^{n-1} f(\tau)\,d\tau,\quad t>0,\, n\in\mathbb{N}.
\end{equation}
This formulation is extended to fractional orders and constitutes the Riemman-Liouville fractional integral expressed for $\alpha\in \mathbb{R}+$ as
\begin{equation}\label{Riemman-Liouville-fractional-integral}
_aD_{t}^{-\alpha}f(t)= \dfrac{1}{\Gamma(\alpha)}\displaystyle\int_{a}^{t} (t-\tau)^{\alpha-1}\,f(\tau)  \,d\tau,\quad t>a.
\end{equation}
where $\Gamma$ is the Euler Gamma function.

A second formulation of fractional differentiation is based on the so-called Hadamard fractional integral as
\begin{equation}\label{Hadamard-fractional-integral}
_aD_{t}^{-\alpha}f(t)=\dfrac{1}{\Gamma(\alpha)}\displaystyle\int_{a}^{t} (\log\dfrac{t}{\tau})^{\alpha-1}\,f(\tau)   \dfrac{d\tau}{\tau}.
\end{equation}
Based on these formulations and analogous the fractional derivatives of arbitrary order $\alpha>0$ becomes a natural requirement. We seek a formulation that remains valid when applied for ordinary integer orders. One formulation is introduced as follows. Let $\alpha>0$ and $m\in\mathbb{N}$ be such that $m-1<\alpha\leq m$. The $\alpha$-derivative of a function $f$ is
\begin{equation}\label{alphaderivative-1}
D^{\alpha}f(t)=D^m J^{m-\alpha}f(t).
\end{equation}
Otherwise,
\begin{equation}\label{alphaderivative-2}
D^\alpha f(t)=
\begin{cases}
\dfrac{d^m}{dt^m}[\dfrac{1}{\Gamma(m-\alpha)}\displaystyle\int_{0}^{t}\dfrac{f(t)}{(t-\tau)^{\alpha+1-m}}],\quad\hbox{for}\;m-1<\alpha<m,\\
\dfrac{d^m}{dt^m} f(t),\quad\hbox{for}\;\alpha=m.
\end{cases}
\end{equation}
The following alternative definition of fractional derivative is originally introduced by Riemann-Liouville and Caputo and formulated using Lagrange's rule for differential operators. Let $\alpha>0$ and $n=[\alpha]$. The $\alpha$-derivative is
\begin{equation}\label{alphaderivative-3}
_aD_{t}^{\alpha}f(t)=\dfrac{d^n}{dt^n}\,_aD_{t}^{-(n-\alpha)}f(t)
\end{equation}
Otherwise, a different alternative has been already formulated by Caputo and states that
\begin{equation}\label{alphaderivative-4}
_a^cD_{t}^{\alpha}f(t)=\dfrac{1}{\Gamma(n-\alpha)} \displaystyle\int_{a}^{t} \dfrac{f^{(n)}(\tau )d\tau}{(t-\tau)^{\alpha+1-n}}
\end{equation}
There are in literature many alternatives and ideas that introduced the fractional derivatives. We refer to the book of Kilbas et al  and Samko et al for more details, calculus and applications \cite{Kilbas}, \cite{Samko}.

In this context, we have the following result.
\begin{lemma}\label{Deriveefractionnaire}
The following assertions hold.
\begin{enumerate}
\item Whenever $p$ and $q$ in $\mathbb{R}_+$ we have
$$
D^q_{\underline{x}}(\underline{b}-\underline{x})= m^q \dfrac{\Gamma(p)}{\Gamma(p-q)}(\underline{b}-\underline{x})^{p-q-1}.
$$	
\item Whenever $r$ and $s$ in $\mathbb{R}_+$ we have
$$
D^r_{\underline{x}}(1+\underline{x})^{s-1}=e^{ir\pi} m^r \dfrac{\Gamma(s)}{\Gamma(s-r)} (1+\underline{x})^{s-r-1}.
$$		
\end{enumerate}
\end{lemma}
\subsection{Fractional Fourier transform}
Let $f$ be in $L^1(\mathbb{R}^m)$. Its Fourier transform denoted usually $\widehat{f}$ or $\mathcal{F}(f)$ is given by the integral transform
$$
\widehat{f}(\eta)=\mathcal{F}(f)(\eta)=\dfrac{1}{(2\pi)^{\frac{m}{2}}}\displaystyle\int_{\mathbb{R}^m}exp(-ix.\eta)f(x)dx,
$$
where $dx$ is the Lebesgues measure on $\mathbb{R}^m$ and $x.\eta$ is the standard inner product of $x$ and $\eta$ in $\mathbb{R}^m$.

The Fractional Fourier Transform has been intensely studied during the last decade. In this section, we will be concerned with the definition and some of its properties. It is based on the well-known Hermite polynomials $\{H_n\}_{n=0}^\infty$ defined on the real line by
$$
H_n(x)=\dfrac{1}{\sqrt{2^nn!\sqrt{\pi}}}(-1)^n\,e^{\frac{x^2}{2}}\dfrac{d^n}{dx^n}e^{-\frac{x^2}{2}}.
$$
It is well-known that these polynomials form an orthonormal basis for $L^2(\mathbb{R}^m)$, and are eigenfunctions of the Fourier transform satisfying $\mathcal{F}H_n=(-i)^nH_n$. Consequently, given $f\in L^2(\mathbb{R}^m)$, its Fourier transform may be written by means of Hermite polynomials as
$$
\mathcal{F}(f)=\displaystyle\sum_{n=0}^{\infty}<f,H_n>(-i)^nH_n
=\displaystyle\sum_{n=0}^{\infty}<f,H_n>e^{-in\frac{\pi}{2}}H_n.
$$
The idea of fractional Fourier transform consists in replacing the fraction $\dfrac{\pi}{2}$ by $\alpha$. Hence, we obtain
$$
\mathcal{F}_af=\displaystyle\sum_{n=0}^{\infty}<f,H_n>e^{-ina}H_n\quad(a\in\mathbb{R}).
$$
Coming back to the Hermite operator $\mathcal{H}$ defined by $\mathcal{H}=\dfrac{1}{2}(-\dfrac{d^2}{dx^2}+x^2-1)$ which has the Hermite polynomials as eigenfunctions, the fractional Fourier transform may be written by means of an exponential operator
$$
\mathcal{F}=e^{-i\frac{\pi}{2}\mathcal{H}}
$$
and the fractional Fourier transform may be similarly be written as
$$
\mathcal{F}_a=e^{-ia\mathcal{H}}.
$$
Backgrounds on fractional Fourier analysis may be found in \cite{Craddock}.
\section{Clifford analysis}
This section is devoted to a brief review on basic concepts on Clifford analysis and fractional calculus as well as the extension of fractional Fourier transforms for the Clifford case. The readers may refer to \cite{Abreuetal}, \cite{Antoine-Murenzi-Vandergheynst}, \cite{Brackx-Schepper-Sommen4}, \cite{Brackx-Schepper-Sommen5}, \cite{Craddock}, \cite{Delanghe}, \cite{DeBie1}, \cite{DeBie2}, \cite{Hitzer}, \cite{Stein-Weiss}.
\subsection{Clifford calculus}
Clifford analysis appeared as a generalization of the complex analysis and Hamiltonians. It extended complex calculus to some type of finite-dimensional associative algebra known as Clifford algebra endowed with suitable operations as well as inner products and norms. It is now applied widely in a variety of fields including geometry and theoretical physics.

Clifford analysis offers a functional theory extending the one of holomorphic functions of complex variable. Starting from the real space $\mathbb{R}^m,\;(m>1)$ (or the complex space $\mathbb{C}^m$) endowed with an orthonormal basis $(e_1,\dots, e_m)$, the Clifford algebra $\mathbb{R}_m$ (or its complexifation $\mathbb{C}_m$) starts by introducing a suitable interior product, let
$$
e_j^2=-1,\quad j=1,\dots,m.
$$
$$
e_je_k+e_ke_j=0,\quad j\neq k,\quad j,k=1,\dots,m.
$$
It is straightforward that is a non-commutative multiplication. Two anti-involutions on the Clifford algebra are important. The conjugation is defined as the anti-involution for which
$$
\overline{e_j}=-e_j,\quad j=1,\dots, m
$$
with the additional rule in the complex case,
$$
\overline{i}=-i.
$$
The inversion is defined as the anti-involution for which
$$
e_j^{+}=e_j,\quad j=1,\dots,m.
$$
This yields a basis of the Clifford algebra ($e_A:A\subset\{1,\dots,m\}$) where $e_{\varnothing}=1$ is the identity element.
As these rules are defined, the euclidean space $\mathbb{R}^m$ is then embedded in the Clifford algebras $\mathbb{R}_m$ and $\mathbb{C}_m$ by identifying the vector $x=(x_1,\dots,x_m)$ with the vector $\underline{x}$ given by
$$
\underline{x}=\displaystyle\sum_{j=1}^{m}e_jx_j.
$$
The product of two vectors is given by
$$
\underline{x}\,\underline{y}=\underline{x}.\underline{y}+\underline{x}\wedge\underline{y}
$$
where
$$
\underline{x}.\underline{y}=-<\underline{x},\underline{y}>=-\displaystyle\sum_{j=1}^{m}x_j\,y_j
$$
and
$$
\underline{x}\wedge\underline{y}=\displaystyle\sum_{j=1}^{m}\displaystyle\sum_{k=j+1}^{m}e_j\,e_k(x_j\,y_k-x_ky_j)
$$
is the wedge product. In particular,
$$
\underline{x}^2=-<\underline{x},\underline{x}>=-|\underline{x}|^2.
$$
\subsection{Clifford monogenic functions}
Let $\Omega$ be an open subset of $\mathbb{R}^m$ or $\mathbb{R}^{m+1}$ and $f:\Omega\rightarrow\mathbb{A}$, where $\mathbb{A}$ is the reel Clifford algebra $\mathbb{R}_{m}$ or its complexification $\mathbb{C}_{m}$. $f$ may be written on the form
$$
f=\displaystyle\sum_{A}f_{A}e_{A}
$$
where the functions $f_A$ are $\mathbb{R}$-valued or $\mathbb{C}$-valued and $(e_A)_A$ is a suitable basis of $\mathbb{A}$. Despite the fact that Clifford analysis generalizes the most important features of classical complex analysis, monogenic functions do not enjoy all properties of holomorphic functions of complex variable.

For instance, due to the non-commutativity of the Clifford algebras, the product of two monogenic functions is in general not monogenic. It is therefore natural to look for specific techniques to construct monogenic functions.

In the literature, there are several techniques available to generate monogenic functions in $\mathbb{R}^{m+1}$ such as the Cauchy-Kowalevski extension (CK-extension) which consists in finding a monogenic extension $g^*$ of an analytic function $g$ defined on a given subset in $\mathbb{R}^{m+1}$ of positive codimension. For analytic functions $g$ on the plane $\{(x_0, \underline{x}) \in\mathbb{R}^{m+1}, \quad x_0 = 0\}$ the problem may be stated as follows:
$$
\left\{\begin{array}{lll}
\mbox{Find}\;\;g^*\in\mathbb{A}\;\;\mbox{such that}\\
\partial_{x_0}g^*=-\partial_{\underline{x}}g^*\;\;in\;\;\mathbb{R}^{m+1},\\
g^*(0, \underline{x}) = g(\underline{x}).
\end{array}
\right.
$$
A formal solution is
\begin{equation}\label{ck-formal-solution}
g^*(x_0,\underline{x})=\exp(-x_0\partial_{\underline{x}}) g(\underline{x})=\displaystyle\sum_{k=0}^{\infty} \dfrac{(-x_0)^k}{k!}
\partial_{\underline{x}}^k g(\underline{x}).
\end{equation}
It may be proved that (\ref{ck-formal-solution}) is a monogenic extension of the function $g$ in $\mathbb{R}^{m+1}$. Moreover, by the uniqueness theorem for monogenic functions this extension is also unique. See \cite{Brackx-Schepper-Sommen0}, \cite{Delanghe}, \cite{Pena}, \cite{Vieira}, \cite{Winkler} and the references therein.

An $\mathbb{R}_m$ or $\mathbb{C}_m$-valued function $F(x_1,\dots,x_m)$, respectively $F(x_0, x_1,\dots,x_m)$ is called right monogenic in an open region of $\mathbb{R}^m$, respectively, or $\mathbb{R}^{m+1}$, if in that region
$$
F\partial_{\underline{x}}=0, \quad respectively\quad F(\partial_{x_0}+\partial_{\underline{x}})=0.
$$
Here $\partial_{\underline{x}}$ is the Dirac operator in $\mathbb{R}^m$:
$$
\partial_{\underline{x}}=\displaystyle\sum_{j=1}^{m} e_j \partial_{x_j},
$$
which splits the Laplacian in $\mathbb{R}^m$
$$
\Delta_m=-\partial_{\underline{x}}^2,
$$
whereas $\partial_{x_0}+\partial_{\underline{x}}$ is the Cauchy-Riemann operator in $\mathbb{R}^{m+1}$ for which
$$
\Delta_{m+1}=(\partial_{x_0}+\partial_{\underline{x}})(\partial_{x_0}+\overline{\partial_{\underline{x}}})
$$
Denoting $S^{m-1}$ the unit sphere in $\mathbb{R}^m$ and introducing spherical co-ordinates in $\mathbb{R}^m$ by
$$
\underline{x}=r\underline{\omega},\quad r=|\underline{x}|\in[0,+\infty[,\;\,\underline{\omega}\in S^{m-1},
$$
the Dirac operator takes the form
$$
\partial_{\underline{x}}=\underline{\omega}\left( \partial_r+\dfrac{1}{r} \Gamma_{\underline{\omega}}\right)
$$
where
$$
\Gamma_{\underline{\omega}}=-\displaystyle\sum_{i<j}e_ie_j(x_i\partial_{x_j}-x_j\partial_{x_i})
$$
is the so-called spherical Dirac operator which depends only on the angular co-ordinates.
\subsection{Clifford Fourier transform}
As for the euclidian case, Fourier analysis is extended to Clifford Fourier analysis \cite{Brackx-Schepper-Sommen2}, \cite{Brackx-Schepper-Sommen5}, \cite{Craddock}. The idea behind the definition of the Clifford Fourier transform originates from the exponential operator representation of the classical Fourier transform by means of Hermite operators. Recall that
$$
\mathcal{F}[f]=\exp(-i\dfrac{\pi}{2}\mathcal{H}_m)[f]=\displaystyle\sum_{n=0}^{\infty}\dfrac{1}{n!}(-i)^n\,\mathcal{H}_m^n[f]	
$$
with $\mathcal{H}_m$ the the classical $m$-dimensional Hermite operator given by
$$
\mathcal{H}_m=\dfrac{1}{2}(\partial_{\underline{x}}^2-\underline{x}^2-m).
$$
To introduce the Clifford analysis character in the Fourier transform, the exponential operator and $\mathcal{H}_m$ have to be replaced by Clifford algebra-valued ones. The starting step consists in factorizing the operator $\mathcal{H}_m$ making use of the factorization of the Laplace operator by the Dirac $\Gamma$ defined by
$$
\Gamma_{\underline{x}}=-\dfrac{1}{2}(\underline{x}\partial_{\underline{x}}-\partial_{\underline{x}}\underline{x}-m).
$$
Next, two Clifford-Hermite operators $H_m^{\pm}$ are introduced,
$$
\mathcal{H}_m^{\pm}=\mathcal{H}_m\pm(\Gamma_{\underline{x}}+\dfrac{m}{2})
$$
\begin{definition}\label{Clifford-Fourier-transform-def}
The Clifford-Fourier transform is the pair of transformations
$$
\mathcal{F}_{m}^+=\exp(-i\dfrac{\pi}{2}\mathcal{H}_m^{+})\quad\hbox{and}\quad\mathcal{F}_{m}^-=\exp(-i\dfrac{\pi}{2}\mathcal{H}_m^-)
$$
\end{definition}
\hskip-17pt Since $\mathcal{H}_m$ commute with $\Gamma$, we have
$$
\mathcal{F}_{m}^+\mathcal{F}_{m}^-=\exp(-i\pi\mathcal{H}_m).
$$
Note that
$$
\mathcal{F}_{m}^\pm=\exp(\pm i(\dfrac{\pi}{2})(\Gamma_{\underline{x}}+\dfrac{m}{2}))\,\mathcal{F}_m
$$
where $\mathcal{F}_m$ is the classical $m$-dimensional Fourier transform which acts by integration against the scalar-valued kernel $K_m(\underline{x},\underline{y})=(2\pi)^{\frac{-m}{2}} \, e^{-i<\underline{x},\underline{y}>}$. As a consequence, we obtain an integral representation for the Clifford-Fourier trasform
$$
\mathcal{F}_{m}^\pm f(\underline{x})=\int_{\mathbb{R}^m} \exp(\pm i\dfrac{\pi}{2})
(\Gamma_{\underline{x}}+\dfrac{m}{2})\, K_{m}(\underline{x},\underline{y})\,f(\underline{y})\,dV(\underline{y})
$$
from which we see that $\mathcal{F}_{m}^\pm$ acts by integration against the Clifford-valued kernel
$$
C_m^{\pm}(\underline{x},\underline{y})=\exp(\pm i\dfrac{\pi}{2})(\Gamma_{\underline{x}}+\dfrac{m}{2})\,K_{m}(\underline{x},\underline{y}).
$$
$dV(\underline{y})$ is the Lebesgue measure on $\mathbb{R}_m$.
\begin{definition}\label{fractional-Clifford-Fourier-transform-def}
The two fractional Clifford-Fourier transform is defined by
$$
\mathcal{F}^{\pm}_{m,a}= \exp(-ia\mathcal{H}_m^{\pm}).
$$
\end{definition}
\hskip-17pt Since we have
$$
\mathcal{F}^{+}_{m,a}\mathcal{F}^{-}_{m,a}=\mathcal{F}_{m,2a},
$$
the fractional Clifford-Fourier transforms act by integration against the kernels
$$
C^{\pm}_{m,a}(\underline{x},\underline{y})= \exp(\pm ia (\Gamma_{\underline{x}}+\dfrac{m}{2})) \, K_{m,a}(\underline{x},\underline{y})
$$
where $K_{m,a} $ is the classical $m$-dimensional fractional Fourier transform.

Throughout this article the Clifford-Fourier transform of $f$ is given by
$$
\mathcal{F}(f(x))(y)=\displaystyle\int_{\mathbb{R}^m} e^{-i<\underline{x},\underline{y}>}\, f(\underline{x}) dV(\underline{x}).
$$
We also denote the inner product of functions in the framework of Clifford analysis by
$$
<f,g>=\displaystyle\int_{\mathbb{R}^m}f(\underline{x})\overline{g(\underline{x})}dV(\underline{x}).
$$
In the present work, we propose to apply such topics to output some generalizations of multidimensional Continuous Wavelet Transform in the context of Clifford analysis.
\section{Generalized Clifford-Jacobi polynomials}
A special case has been developed in \cite{Brackx-Schepper-Sommen4} where the authors tried to point out some wavelets analysis in homogenous euclidean spaces. The main idea applied there was by considering a clifford Heaviside function to decompose the weight $(1-x)^\alpha (1+x)^\beta$ into semi-radial ones and thus to apply then the radial bases of euclidiean spaces by assuming that $\alpha\pm\beta=1$. The basic idea applied there is inspired from there and resumed in the fact that for this case the weight function may be decomposed into a sum of tow parts: one part is the well known radial classical weight $(1+|x|^2)^\alpha$ and a second semi-radial one which leads with the well known rule of derivation of integrals with parameters to the first part. This was the main difference and motivation with our case where such a decomposition is not possible.

This was the main motivation and the main difference with the present paper We instead come back to fractional calculus to overcome the problem of the non radially symmetric weight applied here and did not restrict to the previous case. We will prove instead that fractional calculus may be a good tool to overcome the difficulties crossed in the new context.

In this section we consider the generalization to the Clifford analysis of the classical Jacobi polynomials. Consider the Clifford algebra-valued weight function
$$
\omega_{\alpha,\beta}(\underline{x})=(1-\underline{x})^\alpha(1+\underline{x})^\beta,
$$
with $\alpha, \beta \in \mathbb{R}$ and the monogenic function
$$
F^{*}(t,\underline{x})=\displaystyle\sum_{\ell=0}^{\infty}\dfrac{t^\ell}{\ell!}G_{\ell,m}^{\alpha,\beta}(\underline{x})
\omega_{\alpha-\ell,\beta-\ell}(\underline{x}).
$$
We will evaluate now $\partial_{\underline{x}}(\omega_{\alpha,\beta}(\underline{x}))$. We have
$$
(1-\underline{x})^\alpha= \displaystyle\sum_{n=0}^{\infty} \displaystyle\frac{(\alpha)_n}{n!}(-1)^n \underline{x}^n
$$
and
$$
(1+\underline{x})^\beta= \displaystyle\sum_{n=0}^{\infty} \displaystyle\frac{(\beta)_n}{n!} \underline{x}^n
$$
where for $s\in\mathbb{R}$,
$$
(s)_n=\dfrac{\Gamma(s+1)}{\Gamma(s-n)}.
$$
Thus, we obtain
$$
(1-\underline{x})^\alpha (1+\underline{x})^\beta=\displaystyle\sum_{n=0}^{\infty} a_n(\alpha,\beta) \underline{x}^n
$$
with
$$
a_n(\alpha,\beta) =\displaystyle\sum_{k=0}^{n} \displaystyle\frac{(\beta)_k}{k!} \displaystyle\frac{(\alpha)_{n-k}}{(n-k)!}(-1)^{n-k}.
$$
Next, we use the following technical lemma.
\begin{lem}\label{lemmaevenodd}
For $n\in\mathbb{N}$, we have
$$
\partial_{\underline{x}}(\underline{x}^n)=\gamma_{n,m} \underline{x}^{n-1}.
$$
where,
$$
\gamma_{n,m}=
\begin{cases}
-n \quad \hbox{if} \quad n\quad \mbox{is \,even}.\\
-(m+n-1) \quad \mbox{if}\quad n \quad\,\mbox{is\, odd}.
\end{cases}
$$
\end{lem}
\hskip-17pt Consequently,
$$
\begin{array}{lll}
\partial_{\underline{x}}(\omega_{\alpha,\beta}(\underline{x}))
&=&\displaystyle\sum_{n=0}^{\infty} a_n(\alpha,\beta) \partial_{\underline{x}}(\underline{x}^n)\\
&=& \displaystyle\sum_{n=0}^{\infty} a_n(\alpha,\beta) \gamma_{n,m}\underline{x}^{n-1}.
\end{array}
$$
Otherwise,
$$
\begin{array}{lll}
\partial_{\underline{x}}(\omega_{\alpha,\beta}(\underline{x}))
&=& -\displaystyle\sum_{p=0}^{\infty} a_{2p} (2p) \underline{x}^{2p-1}-
\displaystyle\sum_{p=0}^{\infty} a_{2p+1} (m+2p+1-1)\underline{x}^{2p}\\
&=&-\displaystyle\sum_{n=0}^{\infty} a_{n+1} (n+1) \underline{x}^{n}-\displaystyle\frac{m-1}{\underline{x}} \displaystyle\sum_{n=0}^{\infty} a_{2n+1} \underline{x}^{2n+1}.
\end{array}
$$
Denote
$$
\Gamma_{\alpha,\beta}(\underline{x})=\displaystyle\sum_{n=0}^{\infty} a_{n+1} (n+1) \underline{x}^{n}.
$$
We now use the next Lemma.
\begin{lem}
$$
(n+1)a_{n+1}(\alpha,\beta)=\beta a_n(\alpha,\beta-1)-\alpha a_n(\alpha-1,\beta).
$$	
\end{lem}
\hskip-17pt Hence, we obtain
$$
\begin{array}{lll}
\Gamma_{\alpha,\beta}
&=&\alpha\displaystyle\sum_{n=0}^{\infty}a_{n}(\alpha-1,\beta)\underline{x}^{n}-\beta
\displaystyle\sum_{n=0}^{\infty} a_{n}(\alpha,\beta-1)\underline{x}^{n}\\
&=&\alpha\omega_{\alpha-1,\beta}(\underline{x})-\beta\omega_{\alpha,\beta-1}	(\underline{x}).
\end{array}
$$
As a result,
$$
\partial_{\underline{x}}(\omega_{\alpha,\beta}(\underline{x}))=\alpha \omega_{\alpha-1,\beta}(\underline{x})-\beta\omega_{\alpha,\beta-1}(\underline{x})-\displaystyle\frac{m-1}{2\underline{x}}[\omega_{\alpha,\beta}(\underline{x})-\omega_{\alpha,\beta}(\underline{-x})].
$$
As previously, we will investigate the fact that $F^*(t,\underline{x})$ is monogenic. Observe that
$$
\partial_t F^*(t,\underline{x})=\displaystyle\sum_{\ell=0}^{\infty}\dfrac{t^\ell}{\ell!} G_{\ell+1,m}^{\alpha,\beta}(\underline{x})\, \omega_{\alpha-\ell-1,\beta-\ell-1}(\underline{x})
$$
and
$$
\partial_{\underline{x}}F^*(t,\underline{x})
=\partial_{\underline{x}}G_{\ell,m}^{\alpha,\beta}(\underline{x})\omega_{\alpha-\ell,\beta-\ell}(\underline{x})
+G_{\ell,m}^{\alpha,\beta}(\underline{x})
\partial_{\underline{x}}\left[\omega_{\alpha-\ell,\beta-\ell}(\underline{x})\right].
$$
On the other hand,
$$
\begin{array}{lll}
\partial_{\underline{x}}(\omega_{\alpha-\ell,\beta-\ell}(\underline{x}))
&=&(\alpha-\ell)\omega_{\alpha-\ell-1,\beta-\ell}(\underline{x})-(\beta-\ell)\omega_{\alpha-\ell,\beta-\ell-1}(\underline{x})\\
&&-\displaystyle\frac{m-1}{2\underline{x}}[\omega_{\alpha-\ell,\beta-\ell}(\underline{x})-\omega_{\alpha-\ell,\beta-\ell}(\underline{-x})].
\end{array}
$$
Hence,
$$
\begin{array}{lll}
\,&\;\;&\partial_{\underline{x}}F^*(t,\underline{x})\\
&=&\partial_{\underline{x}}G_{\ell,m}^{\alpha,\beta}(\underline{x})\omega_{\alpha-\ell,\beta-\ell}(\underline{x})\\
&&+G_{\ell,m}^{\alpha,\beta}(\underline{x})\left[(\alpha-\ell) \omega_{\alpha-\ell-1,\beta-\ell}(\underline{x})-(\beta-\ell)\omega_{\alpha-\ell,\beta-\ell-1}(\underline{x})\right.\\
&&\left.-\displaystyle\frac{m-1}{2\underline{x}}[\omega_{\alpha-\ell,\beta-\ell}(\underline{x})-\omega_{\alpha-\ell,\beta-\ell}(\underline{-x})]\right].
\end{array}
$$
Then, writing
$$
\begin{array}{lll}
\,&\;\;&(\partial_{t}+\partial_{\underline{x}})F^{*}(t,\underline{x})\\
&=&G_{\ell+1,m}^{\alpha,\beta}(\underline{x})\omega_{\alpha-\ell-1,\beta-\ell-1}
+\omega_{\alpha-\ell,\beta-\ell}(\underline{x})\partial_{\underline{x}}G_{\ell,m}^{\alpha,\beta}(\underline{x})\\
&+&G_{\ell,m}^{\alpha,\beta}(\underline{x})\left[
(\alpha-\ell) \omega_{\alpha-\ell-1,\beta-\ell}(\underline{x})-(\beta-\ell)\omega_{\alpha-\ell,\beta-\ell-1}(\underline{x})\right.\\
&&\left.-\displaystyle\frac{m-1}{2\underline{x}}[\omega_{\alpha-\ell,\beta-\ell}(\underline{x})-\omega_{\alpha-\ell,\beta-\ell}(\underline{-x})]\right]\\
&=&0,
\end{array}
$$
we get
\begin{equation}\label{5}
G_{\ell+1,m}^{\alpha,\beta}(\underline{x})=-\left[\alpha-\beta+(\alpha+\beta-2\ell)\underline{x}\right]G_{\ell,m}^{\alpha,\beta}(\underline{x})
-(1-\underline{x}^2)\partial_{\underline{x}}G_{\ell,m}^{\alpha,\beta}(\underline{x}).
\end{equation}
This allows to compute $G_{\ell,m}^{\alpha,\beta}(\underline{x})$ recursively. Starting from $G_{0,m}^{\alpha,\beta}(\underline{x})=1$, we deduce that
$$
G_{1,m}^{\alpha,\beta}(\underline{x})=-\left[(\alpha-\beta)+(\alpha+\beta)\underline{x}\right].
$$	
Next, For $\ell=1$,
$$
\begin{array}{lll}
G_{2,m}^{\alpha,\beta}(\underline{x})&=&
\left[(\alpha+\beta)(\alpha+\beta-2+m)\right]\underline{x}^2+
\left[ (\alpha-\beta)(2\alpha+2\beta-2)\right]\underline{x} \\
&&+ (\alpha-\beta)^2-m(\alpha+\beta).	
\end{array}
$$
For $\ell=2$, we get
$$
\begin{array}{lll}
G_{3,m}^{\alpha,\beta}(\underline{x})&=&-[(\alpha+\beta)(\alpha+\beta-2+m)+(\alpha+\beta-2)]\underline{x}^3\\
&&-(\alpha-\beta)[(\alpha+\beta)(\alpha+\beta-2+m)+(2\alpha+2\beta-2)(\alpha+\beta-4-m)]\underline{x}^2\\
&&-[(\alpha-\beta)^2(2\alpha+2\beta-2)
+(\alpha+\beta)[(\alpha-\beta)^2-m(\alpha+\beta)-2(\alpha+\beta-2+m)]\underline{x}\\
&&-(\alpha-\beta)[(\alpha+\beta)^2-m(\alpha+\beta)-m
(2\alpha+2\beta-2)].
\end{array}
$$
\begin{proposition}
The generalized Clifford-Jacobi polynomials can be written as
$$
G_{\ell,m}^{\alpha,\beta}(\underline{x})=
(-1)^{\ell}\omega_{\ell-\alpha,\ell-\beta}(\underline{x})\partial_{\underline{x}}^\ell\omega_{\alpha,\beta}(\underline{x}).
$$
\end{proposition}
\hskip-17pt\textbf{Proof.} The proof is based on the recurrencxe principle. We will sketch it just for the first two steps. The rest is left to the reader. For $\ell=1$, we have
$$
\begin{array}{lll}
\partial_{\underline{x}}\omega_{\alpha,\beta}(\underline{x})
&=&\alpha\,\omega_{\alpha-1,\beta}(\underline{x})-\beta\omega_{\alpha,\beta-1}(\underline{x})-\displaystyle\frac{m-1}{2\underline{x}}[\omega_{\alpha,\beta}(\underline{x})-\omega_{\alpha,\beta}(-\underline{x})]\\
&=&(-1)\omega_{\alpha-1,\beta-1}(\underline{x})\left[(\alpha-\beta)+(\alpha+\beta)\underline{x}\right]\\
&=&(-1)\omega_{\alpha-1,\beta-1}(\underline{x}) G_{1,m}^{\alpha,\beta}(\underline{x}).
\end{array}
$$
Which means that
$$
G_{1,m}^{\alpha,\beta}=(-1)\,\omega_{1-\alpha,1-\beta}(\underline{x}) \partial_{\underline{x}}\,\omega_{\alpha,\beta}(\underline{x}).
$$
For $\ell=2$, we shall get
\begin{equation}\label{g2}
G_{2,m}^{\alpha,\beta}= (-1)^2\omega_{2-\alpha,2-\beta}(\underline{x})\partial_{\underline{x}}^{(2)}\omega_{\alpha,\beta}(\underline{x}).
\end{equation}
Indeed,
$$
\begin{array}{lll}
\partial_{\underline{x}}^{(2)}\omega_{\alpha,\beta}(\underline{x})
&=&\partial_{\underline{x}}(\alpha\,\omega_{\alpha-1,\beta}(\underline{x})-\beta\omega_{\alpha,\beta-1}(\underline{x})\\
&&-\displaystyle\frac{m-1}{2\underline{x}}[\omega_{\alpha,\beta}(\underline{x})-\omega_{\alpha,\beta}(-\underline{x})])\\
&=&\partial_{\underline{x}}\left(\alpha\,\omega_{\alpha-1,\beta}(\underline{x})-\beta\omega_{\alpha,\beta-1}(\underline{x})-(m-1)\displaystyle\sum_{n=0}^{\infty}a_{2n+1} \underline{x}^{2n}\right)\\
&=& \alpha\,\partial_{\underline{x}}(\omega_{\alpha-1,\beta}(\underline{x}))-\beta \partial_{\underline{x}}(\omega_{\alpha,\beta-1}(\underline{x}))\\
&&-(m-1)\partial_{\underline{x}}\left(\displaystyle\sum_{n=0}^{\infty}a_{2n+1} \underline{x}^{2n}\right)\\
&=&\alpha[(\alpha-1)\omega_{\alpha-2,\beta}(\underline{x})-\beta\omega_{\alpha-1,\beta-1}(\underline{x})\\
&-&\displaystyle\frac{m-1}{2\underline{x}}(\omega_{\alpha-1,\beta}(\underline{x})-\omega_{\alpha-1,\beta}(-\underline{x}))]\\
&-&\beta[\alpha\omega_{\alpha-1,\beta-1}(\underline{x})-(\beta-1)\omega_{\alpha,\beta-2}(\underline{x})\\
&-&\displaystyle\frac{m-1}{2\underline{x}}(\omega_{\alpha,\beta-1}(\underline{x})-\omega_{\alpha,\beta-1}(-\underline{x}))]\\
&+&\displaystyle\frac{m-1}{\underline{x}^2}(\omega_{\alpha,\beta}(\underline{x})-\omega_{\alpha,\beta}(\underline{-x}))]\\
&=&(-1)^2\omega_{\alpha-2,\beta-2}(\underline{x})[\alpha(\alpha-1)(1+\underline{x})^2-\alpha\beta(1-\underline{x}^2)\\
&-&\alpha\displaystyle\frac{m-1}{2\underline{x}}[(1-\underline{x})(1+\underline{x})^2-(1-\underline{x})^2(1+\underline{x})]\\
&-&\alpha\beta(1-\underline{x}^2)+\beta(\beta-1)(1-\underline{x})^2\\
&+&\beta\displaystyle\frac{m-1}{2\underline{x}}[(1-\underline{x})^2(1+\underline{x})-(1-\underline{x})(1+\underline{x})^2]\\
&+&\displaystyle\frac{m-1}{\underline{x}^2}[(1-\underline{x})^2(1+\underline{x})^2-(1-\underline{x})^2(1+\underline{x})^2]]\\
&=&(-1)^2\omega_{\alpha-2,\beta-2}(\underline{x})[\alpha(\alpha-1)+2\alpha(\alpha-1)\underline{x}\\
&+&\alpha(\alpha-1)\underline{x}^2-2\alpha\beta+2\alpha\beta\underline{x}^2-\alpha(m-1)\\
&+&\alpha(m-1)\underline{x}^2+\beta(\beta-1)-2\beta(\beta-1)\underline{x}\\
&+&\beta(\beta-1)\underline{x}^2+\beta(m-1)\underline{x}^2-\beta(m-1)].
\end{array}
$$
By regrouping the quantities with $\underline{x}^2,\underline{x}$, we get
$$
\partial_{\underline{x}}^{(2)}\omega_{\alpha,\beta}(\underline{x})=(-1)^2 \omega_{\alpha-2,\beta-2}(\underline{x})
[G_{2,m}^{\alpha,\beta}(\underline{x})]
$$
where $G_{2,m}^{\alpha,\beta}(\underline{x})$ is as in (\ref{g2}).
\begin{proposition} Let the  integral
$$
I_{\ell,t,p}^{\alpha,\beta}=\displaystyle\int_{\mathbb{R}^m}\underline{x}^\ell G_{t,m}^{\alpha+p,\beta+p}(\underline{x})\,\omega_{\alpha, \beta}(\underline{x})\,dV(\underline{x}).
$$
Then, whenever $2t<1-\alpha-\beta-m$ and  $\ell<t$, we have the orthogonality relation
$$
I_{\ell,t,t}^{\alpha,\beta}=0.
$$
\end{proposition}
\hskip-17pt\textbf{Proof.} Denote,
$$
I_{\ell,t}=\displaystyle\int_{\mathbb{R}^m}\underline{x}^\ell\partial_{\underline{x}}^t\omega_{\alpha+t,\beta+t}(\underline{x}) dV(\underline{x}).
$$
Using Stokes's theorem, we obtain
$$
\begin{array}{lll}
\;&\;&\displaystyle\int_{\mathbb{R}^m}\underline{x}^\ell\,G_{t,m}^{\alpha+t,\beta+t}(\underline{x})\omega_{\alpha,\beta}(\underline{x})dV(\underline{x})\\
&=&(-1)^t\displaystyle\int_{\mathbb{R}^m}\underline{x}^\ell\omega_{t-\alpha-t,t-\beta-t}(\underline{x})
[\partial_{\underline{x}}^t\omega_{\alpha+t,\beta+t}(\underline{x})]\omega_{\alpha,\beta}(\underline{x}) dV(\underline{x})\\
&=&(-1)^t\displaystyle\int_{\mathbb{R}^m}\underline{x}^\ell\partial_{\underline{x}}^t\omega_{\alpha+t,\beta+t}(\underline{x}) dV(\underline{x})\\
&=&(-1)^t\displaystyle\int_{\mathbb{R}^m}\underline{x}^\ell\partial_{\underline{x}}\left[\partial_{\underline{x}}^{t-1}\omega_{\alpha+t,\beta+t}(\underline{x})\right] dV(\underline{x})\\
&=&(-1)^t\left[\displaystyle\int_{\partial\mathbb{R}^m}\underline{x}^\ell
\partial_{\underline{x}}^{t-1}\omega_{\alpha+t,\beta+t}(\underline{x})\partial\Gamma(\underline{x})
-\displaystyle\int_{\mathbb{R}^m}\partial_{\underline{x}}(\underline{x}^\ell) \partial_{\underline{x}}^{t-1}\omega_{\alpha+t,\beta+t}(\underline{x}) dV(\underline{x})\right].
\end{array}
$$
Denote now
$$
I=\displaystyle\int_{\partial\mathbb{R}^m}\underline{x}^\ell\partial_{\underline{x}}^{t-1}\omega_{\alpha+t,\beta+t}(\underline{x})\partial\Gamma(\underline{x})
$$
and
$$
II=\displaystyle\int_{\mathbb{R}^m}\partial_{\underline{x}}(\underline{x}^\ell) \partial_{\underline{x}}^{t-1}\omega_{\alpha+t,\beta+t}(\underline{x}) dV(\underline{x}).
$$
The integral $I$ vanishes due to the assumption $2t<1-\alpha-\beta-m $.
Now, using Lemma \ref{lemmaevenodd}, we evaluate $II$. Indeed,
$$
\begin{array}{lll}
II&=&\gamma_{l,m}\displaystyle\int_{\mathbb{R}^m}\underline{x}^{\ell-1}
\partial_{\underline{x}}^{t-1}\,\omega_{\alpha+t,\beta+t}(\underline{x})\,dV(\underline{x})\\
&=&\gamma_{l,m}I_{\ell-1,t-1}.
\end{array}
$$
Hence we obtain
$$
\begin{array}{lll}
\displaystyle\int_{\mathbb{R}^m}\underline{x}^\ell G_{t,m}^{\alpha+t,\beta+t}(\underline{x})\omega_{\alpha,\beta}(\underline{x})dV(\underline{x})
&=&(-1)^{t+1}\gamma_{l,m}I_{\ell-1,t-1}\\
&=&(-1)^{t+1}\gamma_{l,m}[(-1)^{t}\gamma_{l-1,m}I_{\ell-2,t-2}]\\
&=&(-1)^{2t+1}\gamma_{l,m}\gamma_{l-1,m}I_{\ell-2,t-2}\\
&\vdots&\\
&=&C(m,\ell,t)I_0\\
&=&0.
\end{array}
$$
The constant $C(m,\ell,t)$ above is given by
$$
C(m,\ell,t)=(-1)^{ml+1}\displaystyle\prod_{k=0}^{m}\gamma_{k,m}.
$$
\section{The Generalized Clifford-Jacobi continuous wavelet transform}
Now we are able to introduce a new class of wavelets relative to the generalized Clifford-Jacobi polynomials. Recall that for $0<t<\dfrac{1-\alpha-\beta-m}{2}$ we have
$$
\displaystyle\int_{\mathbb{R}^m} G_{t,m}^{\alpha+t,\beta+t}(\underline{x})\omega_{\alpha,\beta}(\underline{x})dV(\underline{x})=0.
$$
Consequently, we get the following definition.
\begin{definition}\label{GeneralizedClifford-JacobiWaveletMother} The generalized Clifford-Jacobi wavelet mother is defined by
$$
\psi_{\ell,m}^{\alpha,\beta}(\underline{x})= G_{\ell,m}^{\alpha+\ell,\beta+\ell}(\underline{x})\omega_{\alpha,\beta}(\underline{x}).
$$
\end{definition}
Furthermore, the wavelet $\psi_{\ell,m}^{\alpha,\beta}(\underline{x})$ have vanishing moments as is shown in the next proposition.
\begin{proposition} The following assertions hold.
\begin{enumerate}
\item The wavelet $\psi_{\ell,m}^{\alpha,\beta}$ have vanishing moments if the condition
$\forall k,\,0\leq k\leq-(m+\ell+\alpha+\beta)$ and $0\leq k\leq\ell$ is fulfilled. More precisely,
\begin{equation}\label{8}
\displaystyle\int_{\mathbb{R}^m}\underline{x}^k\psi_{\ell,m}^{\alpha,\beta}(\underline{x})dV(\underline{x})=0,\quad\hbox{for}\quad
0\leq\,k\leq-(m+\ell+\alpha+\beta)\quad\hbox{and}\quad 0\leq k\leq\ell.
\end{equation}
\item The Clifford-Fourier transform of $\psi_{\ell,m}^{\alpha,\beta}$ is given by
\begin{equation}\label{9}
\widehat{\psi_{\ell,m}^{\alpha,\beta}}(\underline{u})=(-i\underline{u})^{\ell} C_{p,q}e^{-ir\pi}C_{r,s}\displaystyle\sum_{k=0}^{N_{\alpha}} \displaystyle\sum_{n=0}^{N_{\beta}}C_{N_{\alpha}}^k C_{N_{\beta}}^n (-1)^k\widetilde{C_{k,q}}\widetilde{C_{n,r}}K_{k,n}^{q,r}(\underline{u})
\end{equation}
where
$$
C_{p,q}=m^{-q} \dfrac{\Gamma(p-q)}{\Gamma(p)},\;\widetilde{C_{k,q}}=\dfrac{k! m^{k-q}}{\Gamma(k+1-q)},\;N_{\alpha}=N_{\alpha,l}^q=\alpha+\ell+q
$$
and
$$
K_{k,n}^{q,r}(\underline{u})=e^{i\dfrac{\pi}{2}(k+n-q-r)}\dfrac{(2\pi)^{\frac{m}{2}}K_m}{\Gamma(m)} \Gamma(k+n-q-r+m-1) |\underline{u}|^{q+r-k-n-m}.
$$
\end{enumerate}
\end{proposition}
\textbf{Proof.}
$$
\begin{array}{lll}
\widehat{\psi_{\ell,m}^{\alpha,\beta}}(\underline{u})&=&\displaystyle\int_{\mathbb{R}^m}\psi_{\ell,m}^{\alpha,\beta}(\underline{x}) e^{-i\underline{x}.\underline{u}}\,d\underline{x}\\
&=&(-1)^\ell\displaystyle\int_{\mathbb{R}^m}\partial_{\underline{x}}^\ell\left(\omega_{\alpha+\ell,\beta+\ell}(\underline{x})\right)
e^{-i\underline{x}.\underline{u}}\,d\underline{x}\\
&=& (-1)^\ell\displaystyle\int_{\mathbb{R}^m}\omega_{\alpha+\ell,\beta+\ell}(\underline{x})
e^{-i\underline{x}.\underline{u}}(i\underline{u})^\ell \,d\underline{x}\\
&=&(-1)^\ell\,(i\underline{u})^\ell\displaystyle\int_{\mathbb{R}^m}\omega_{\alpha+\ell,\beta+\ell}(\underline{x})\, e^{-i\underline{x}.\underline{u}}\,d\underline{x}\\
&=&(-1)^{\ell}\,(i\underline{u})^{\ell}\displaystyle\int_{\mathbb{R}^m}
(1-\underline{x})^{\alpha+\ell}(1+\underline{x})^{\beta+\ell}\,e^{-i\underline{x}.\underline{u}} \,d\underline{x}\\
&=&(-1)^\ell (i\underline{u})^\ell \widehat{\omega_{\alpha+\ell,\beta+\ell}}(\underline{u}).
\end{array}
$$
Applying Lemma \ref{Deriveefractionnaire}, for $b=1,\, p-q-1=\alpha+\ell$ we get
$$
(D^q_{\underline{x}}(1-\underline{x})^{p-1})= m^q \dfrac{\Gamma(p)}{\Gamma(p-q)}(1-\underline{x})^{\alpha+\ell}.
$$
Therefore,
$$
(1-\underline{x})^{\alpha+\ell}=C_{p,q}(D^q_{\underline{x}}(1-\underline{x})^{p-1}).
$$
Consequently,
$$
\begin{array}{lll}
\widehat{\psi_{\ell,m}^{\alpha,\beta}}(\underline{u})&=&(-i\underline{u})^{\ell}C_{p,q}\displaystyle\int_{\mathbb{R}^m}D_{\underline{x}}^q(1-\underline{x})^{p-1} (1+\underline{x})^{\beta+\ell} e^{-i\underline{x}.\underline{u}} d\underline{x}\\
&=&  (-i\underline{u})^{\ell}C_{p,q} \displaystyle\int_{\mathbb{R}^m}D_{\underline{x}}^q(1-\underline{x})^{\alpha+\ell+q} (1+\underline{x})^{\beta+\ell} e^{-i\underline{x}.\underline{u}} d\underline{x}.
\end{array}
$$
For $q=[\alpha]+1-\alpha$, it holds that $\alpha+\ell+q \in \mathbb{N}$. Hence,
$$
\begin{array}{lll}
\widehat{\psi_{\ell,m}^{\alpha,\beta}}(\underline{u})&=&
(-i\underline{u})^{\ell}C_{p,q}\displaystyle\sum_{k=0}^{\alpha+\ell+q} C_{\alpha+\ell+q}^k(-1)^k\displaystyle\int_{\mathbb{R}^m}D_{\underline{x}}^q \underline{x}^k(1+\underline{x})^{\beta+\ell}e^{-i\underline{x}.\underline{u}} d\underline{x}\\
&=& (-i\underline{u})^{\ell}C_{p,q} \displaystyle\sum_{k=0}^{\alpha+\ell+q} C_{\alpha+\ell+q}^k (-1)^kH_{k,q}\displaystyle\int_{\mathbb{R}^m}\underline{x}^{k-q} (1+\underline{x})^{\beta+\ell} e^{-i\underline{x}.\underline{u}} d\underline{x}\\
&=& (-i\underline{u})^{\ell}C_{p,q} \displaystyle\sum_{k=0}^{\alpha+\ell+q} C_{\alpha+\ell+q}^k (-1)^k H_{k,q}I_{q,\beta}^{l,k}(\underline{u}),
\end{array}
$$
because of the fact that
$$
D_{\underline{x}}^q (\underline{x}^k)=H_{k,q}x^{k-q}
$$
where
$$
H_{k,q}=\prod_{r=0}^{q-1} \gamma_{k-r,m}
$$
and where we set
$$
I_{q,\beta}^{l,k}(\underline{u})=\displaystyle\int_{\mathbb{R}^m}\underline{x}^{k-q}(1+\underline{x})^{\beta+\ell}e^{-i\underline{x}.\underline{u}} d\underline{x}.
$$
We now evaluate this last integral. To do it, we need again Lemma \ref{Deriveefractionnaire}. For $s-r-1=\beta+\ell$, we get
$$
(1+\underline{x})^{\beta+\ell}=e^{-ir\pi}C_{r,s}D_{\underline{x}}^r (1+\underline{x})^{s-1}
$$
where $s-1=\beta+\ell+r$. For $r=[\beta]+1-\beta$ we get $s-1=\ell+[\beta]+1\in \mathbb{N}$. Consequently,
$$
\begin{array}{lll}
(1+\underline{x})^{\beta+\ell}
&=& e^{-ir\pi}C_{r,s} D_{\underline{x}}^r \displaystyle\sum_{n=0}^{\beta+\ell+r}C_{\beta+\ell+r}^n\,\underline{x}^n\\
&=& e^{-ir\pi}C_{r,s}\displaystyle\sum_{n=0}^{\beta+\ell+r}C_{\beta+\ell+r}^n D_{\underline{x}}^r \,\underline{x}^n\\
&=& e^{-ir\pi}C_{r,s} \displaystyle\sum_{n=0}^{\beta+\ell+r}C_{\beta+\ell+r}^n H_{n,r}\,\underline{x}^{n-r}.
\end{array}
$$
As a result,
$$
I_{q,\beta}^{l,k}(\underline{u})=e^{-ir\pi}C_{r,s} \displaystyle\sum_{n=0}^{\beta+\ell+r}C_{\beta+\ell+r}^n H_{n,r} \displaystyle\int_{\mathbb{R}^m}\underline{x}^{k+n-q-r} e^{-i\underline{x}\underline{u}} d\underline{x}.
$$
Denote now
$$
K_{k,n}^{q,r}(\underline{u})=\displaystyle\int_{\mathbb{R}^m}\underline{x}^{k+n-q-r} e^{-i\underline{x}\underline{u}} d\underline{x}.
$$
It is straightforward that
$$
K_{k,n}^{q,r}(\underline{u})=(-i)^{q+r-k-n}D_{\underline{u}}^{k+n-q-r}\displaystyle\int_{\mathbb{R}^m}e^{-i\underline{x}\underline{u}}d\underline{x}.
$$
Therefore,
$$
\widehat{\psi_{\ell,m}^{\alpha,\beta}}(\underline{u})=(-i\underline{u})^{\ell}C_{p,q}e^{-ir\pi}C_{r,s}\displaystyle\sum_{k=0}^{N_{\alpha}} \displaystyle\sum_{n=0}^{N_{\beta}}C_{N_{\alpha}}^k C_{N_{\beta}}^n (-1)^k H_{k,q} H_{n,r} K_{k,n}^{q,r}(\underline{u}).
$$
We now evaluate the quantity $K_{k,n}^{q,r}(\underline{u})$. So denote
$$
K_{k,n}^{q,r}(\underline{u})=(-i)^{q+r-k-n}D_{\underline{u}}^{k+n-q-r}\left(\displaystyle\int_{\mathbb{R}^m} e^{-i\underline{x}\underline{u}} d\underline{x}\right).
$$
and
$$
J(\underline{u})=\displaystyle\int_{\mathbb{R}^m} e^{-i\underline{x}\underline{u}} d\underline{x}.
$$
We have the following technical lemma.
\begin{lemma}
$$
\displaystyle\int_{S^{m-1}}e^{-i<r\underline{\omega},\rho\underline{\xi}>}d\underline{\xi}=\dfrac{(2\pi)^{\frac{m}{2}}J_{\frac{m}{2}-1}
(r\rho)}{(r\rho)^{\frac{m}{2}-1}}.
$$
\end{lemma}
\hskip-17pt By applying the spherical coordinates
$$
\underline{x}= t\underline{\omega}, \quad \underline{u}=\rho \underline{\xi}, \quad t=|\underline{x}|, \quad \rho=|\underline{u}|, \quad \omega,\xi \in S^{m-1}
$$
we get
$$
\begin{array}{lll}
J(\underline{u})&=&\displaystyle\int_{0}^{+\infty}t^{m-1}dt\displaystyle\int_{S^{m-1}}e^{-it\rho\underline{\omega}\underline{\xi}}d\tau(\underline{\omega})\\
&=&\displaystyle\int_{0}^{+\infty}t^{m-1}(2\pi)^{\frac{m}{2}}\dfrac{J_{\frac{m}{2}-1}(t\rho)}{t^{\frac{m}{2}-1}\rho^{\frac{m}{2}-1}}dt\\
&=&(2\pi)^{\frac{m}{2}}\displaystyle\int_{0}^{+\infty}(\dfrac{t}{\rho})^{\frac{m}{2}-1}J_{\frac{m}{2}-1}(t\rho)tdt\\
&=&(2\pi)^{\frac{m}{2}}\rho^{-m}\displaystyle\int_{0}^{+\infty}x^{\frac{m}{2}}J_{\frac{m}{2}-1}(x)dx.
\end{array}
$$
Denoting now
$$
K_m=\displaystyle\int_{0}^{+\infty}x^{\frac{m}{2}}J_{\frac{m}{2}-1}(x)dx
$$
we obtain
$$
K_{k,n}^{q,r}(\underline{u})=(-i)^{q+r-k-n}(2\pi)^{\frac{m}{2}}K_mD_{\rho}^{k+n-q-r}(\rho^{-m}).
$$
Now, as
$$
D^{\nu}_{\rho}(\rho^{-m})= \dfrac{\Gamma(m+\nu-1)}{\Gamma(m)} e^{i\nu\pi} \rho^{-m-\nu},
$$
we get
$$
K_{k,n}^{q,r}(\underline{u})=\Gamma_{q,r}^{m,n}(k)|\underline{u}|^{q+r-k-n-m},
$$
where
$$
\Gamma_{q,r}^{m,n}(k)=e^{i\dfrac{\pi}{2}(k+n-q-r)}\dfrac{(2\pi)^{\frac{m}{2}}K_m}{\Gamma(m)} \Gamma(k+n-q-r+m-1) .
$$
Now, we introduce the generalized Clifford-Jacobi continuous wavelet transform.
\begin{definition}
For $a>0$ and $\underline{b}\in\mathbb{R}_m$, the $(a,\underline{b})$-copy of the wavelet mother $\psi_{\ell,m}^{\alpha,\beta}$ is defined by
\begin{equation} \label{11}
^a_{\underline{b}}{\psi}_{\ell,m}^{\alpha,\beta}(\underline{x})=a^{-\frac{m}{2}}{\psi}_{\ell,m}^{\alpha,\beta}(\dfrac{\underline{x}-\underline{b}}{a}).
\end{equation}
The generalized Clifford-Jacobi CWT of a function $f\in L_2(\mathbb{R}^m)$ is defined by
$$
C_{a,\underline{b}}(f)=<^a_{\underline{b}}\!\!{\psi}_{\ell,m}^{\alpha,\beta},f>.
$$
\end{definition}
The following theorem guaranties the construction of the analyzed function $f\in L_2(\mathbb{R}^m)$ from its wavelet transform.
\begin{theorem}\label{TheoremReconstruction} Let ${\psi}_{\ell,m}^{\alpha,\beta}$ be an analyzing wavelet as in Definition \ref{GeneralizedClifford-JacobiWaveletMother}. The following assertions hold.
\begin{enumerate}
\item ${\psi}_{\ell,m}^{\alpha,\beta}$ is admissible in the sense that
$$
{\mathcal{A}_{\ell,m}^{\alpha,\beta}}=\dfrac{1}{\omega_m}\displaystyle\int_{\mathbb{R}^m}\left|\widehat{\psi_{\ell,m}^{\alpha,\beta}}(\underline{x})\right|^2 \dfrac{dV(\underline{x})}{|\underline{x}|^m}<+\infty,
$$
where $\omega_m$ is the area of the unit sphere $S^{m-1}$ in $\mathbb{R}^m$.
\item The analyzed function $ff\in L_2(\mathbb{R}^m)$ may be reconstructed in the $L_2$-sense by
$$
f(x)=\dfrac{1}{{\mathcal{A}_{\ell,m}^{\alpha,\beta}}}\displaystyle\int_{a>0}\displaystyle\int_{b\in\mathbb{R}^m}C_{a,\underline{b}}(f)
\psi_{\ell,m}^{\alpha,\beta}\left(\dfrac{\underline{x}-\underline{b}}{a}\right)\dfrac{da\,dV(\underline{b})}{a^{m+1}}.
$$
\end{enumerate}
\end{theorem}
The proof reposes on the following result.
\begin{lemma}
Consider the inner product
$$
<C_{a,\underline{b}}(f),C_{a,\underline{b}}(g)>
=\dfrac{1}{\mathcal{A}_{\ell,m}^{\alpha,\beta}}\displaystyle\int_{\mathbb{R}^m}\displaystyle\int_{0}^{+\infty}C_{a,\underline{b}}(f)\overline{C_{a,\underline{b}}(g)}\dfrac{da}{a^{m+1}}dV(\underline{b}).
$$
Then, the following Parseval formula hold
$$
<C_{a,\underline{b}}(f),C_{a,\underline{b}}(g)>=<f,g>.
$$
\end{lemma}
\hskip-17pt\textbf{Proof of Theorem \ref{TheoremReconstruction}.} Assertion \textbf{1.} is based on the asymptotic behaviour of Bessel functions. (Analogue result is already checked in \cite{Brackx-Schepper-Sommen1}.\\
\textbf{2.} Using the Clifford Fourier transform we observe that
$$
C_{a,\underline{b}}(f)(\underline{b})=\widetilde{a^{\frac{m}{2}}\widehat{\widehat{f}(\underline{.})\widehat{\psi_{\ell,m}^{\alpha,\beta}}(a\underline{.})}}(\underline{b}),
$$
where, $\widetilde{h}(\underline{u})=h(-\underline{u})$, $\forall\,h$. Thus,
$$
C_{a,\underline{b}}(f)\overline{C_{a,\underline{b}}(g)}= \widehat{\left(\widehat{f}(\underline{.})a^{\frac{m}{2}} \widehat{\psi_{\ell,m}^{\alpha,\beta}}(a\underline{.})\right)}(-\underline{b})
\overline{\widehat{\left(\widehat{g}(\underline{.})a^{\frac{m}{2}} \widehat{\psi_{\ell,m}^{\alpha,\beta}}(a\underline{.})\right)}}(-\underline{b}).
$$
Consequently,
$$
\begin{array}{lll}
<C_{a,\underline{b}}(f),C_{a,\underline{b}}(g)>
&=&\displaystyle\int_{a>0}\displaystyle\int_{\mathbb{R}^m}\widehat{\widehat{f}(\underline{.})a^{\frac{m}{2}}\widehat{\psi_{\ell,m}^{\alpha,\beta}}(a\underline{.})}\,\overline{\widehat{\widehat{g}(\underline{.})a^{\frac{m}{2}}
\widehat{\psi_{\ell,m}^{\alpha,\beta}}(a\underline{.})}}\dfrac{da\,dV(\underline{b})}{a^{m+1}}\\
&=&\displaystyle\int_{a>0}\displaystyle\int_{\mathbb{R}^m}\widehat{f}(\underline{b})\overline{\widehat{g}(\underline{b})}\dfrac{a^m|
\widehat{\psi_{\ell,m}^{\alpha,\beta}}(a\underline{b})|^2}{a^{m+1}}\,da\,dV(\underline{b})\\
&=&{\mathcal{A}_{\ell,m}^{\alpha,\beta}}\displaystyle\int_{\omega\in\mathbb{R}^m}\widehat{g}(b)\overline{\widehat{g(\underline{b})}}dV(\underline{b})\\
&=&{\mathcal{A}_{\ell,m}^{\alpha,\beta}}<\widehat{f},\widehat{g}>\\
&=&<f,g>.
\end{array}
$$
\section{Conclusion}
In the present work a large class of Clifford-Jacobi polynomials has been developed relatively to general Jacobi weights. Next, new associated wavelets in the Clifford context have been introduced using fractional calculus such as derivatives and fractional Fourier transform on Clifford algebras. The new classes of polynomials as well as wavelets extend the works of Brackx and his collaborators in \cite{Brackx-Schepper-Sommen1}, \cite{Brackx-Schepper-Sommen2}, \cite{Brackx-Schepper-Sommen3}, \cite{Brackx-Schepper-Sommen4} and \cite{Brackx-Schepper-Sommen5}.

\end{document}